\documentclass[12pt]{amsart}

\usepackage{graphicx}
\usepackage[latin1]{inputenc}
\usepackage{xcolor}

\newtheorem{theorem}{Theorem}[section]
\newtheorem{corollary}[theorem]{Corollary}
\newtheorem{lemma}[theorem]{Lemma}

\newtheorem{definition}[theorem]{Definition}
\newtheorem{remark}[theorem]{Remark}
\newtheorem{proposition}[theorem]{Proposition}
\newtheorem{example}[theorem]{Example}

\def\R{{\mathbb R}}

\begin{document}

\title{The groups $G_{n}^{k}$, $2n$-gon tilings, and stacking of cubes}

\author{Seongjeong Kim}

\address{Jilin university, Changchun, China \\
kimseongjeong@jlu.edu.cn}

\author{Vassily Olegovich Manturov}

\address{Moscow Institute of Physics and Technology, Moscow 140180, Russia \\
vomanturov@yandex.ru}

\maketitle

\begin{abstract}
In the present paper we discuss three ways of looking at rhombile tilings: stacking 3-dimensional cubes, elements of groups, and configurations of lines and points.

\end{abstract}

\vspace{0.40in} {\em Keywords:\/} cube, hexagon, octagon, flip, line configuration, projective plane, duality, desargues, configuration space, braid, 3-manifold, spine, rhombile tiling, surface tiling, octagon relation, $G_{n}^{3}$ group. 

\section{Introduction}

In mathematics rhombile tilings of polygons and flip operations on them
were widely studied in \cite{HenriquesSpeyer,DThurston}.
 
\begin{figure}[h!]
\begin{center}
 \includegraphics[width =4cm]{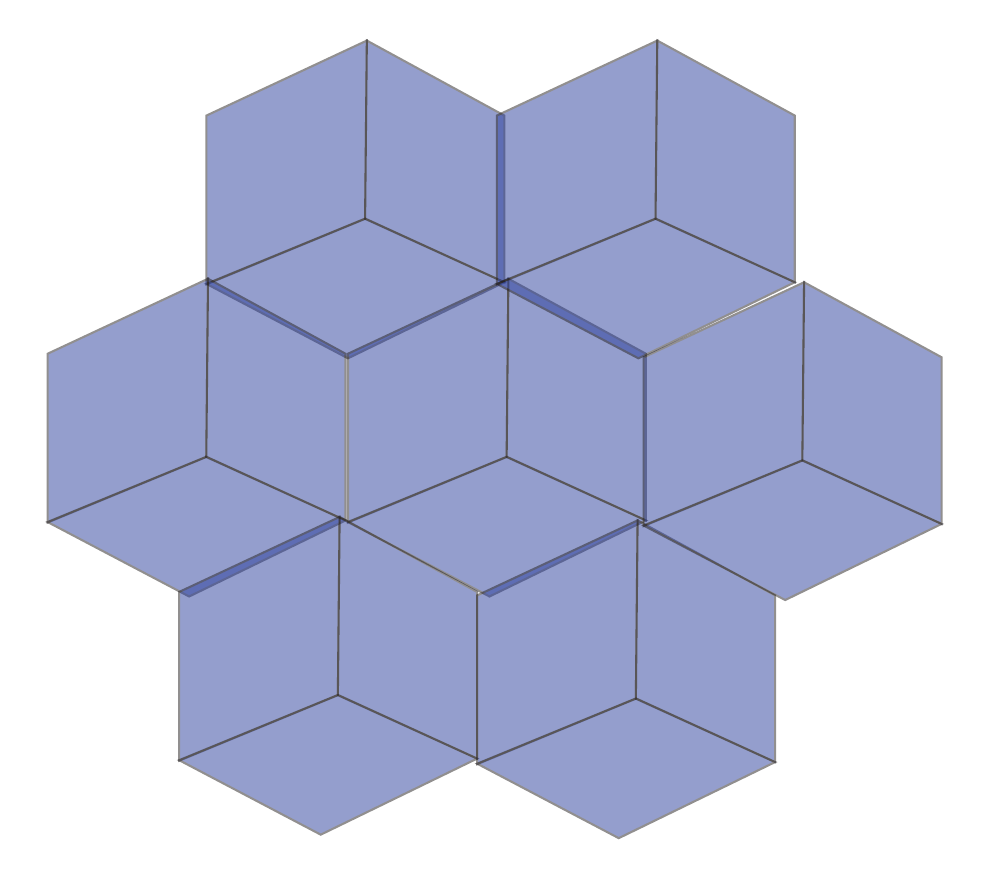}

\end{center}
 \caption{Cubes in 3-space or rhombile tiling of a plane}\label{fig:cubes-tiling}
\end{figure}

A beautiful picture Fig. \ref{fig:cubes-tiling} can be treated at least in two ways:
one can see either cubes in 3-space or a tiling of a planar figure into three types
of rhombi. 

\begin{figure}[h!]
\begin{center}
 \includegraphics[width =10cm]{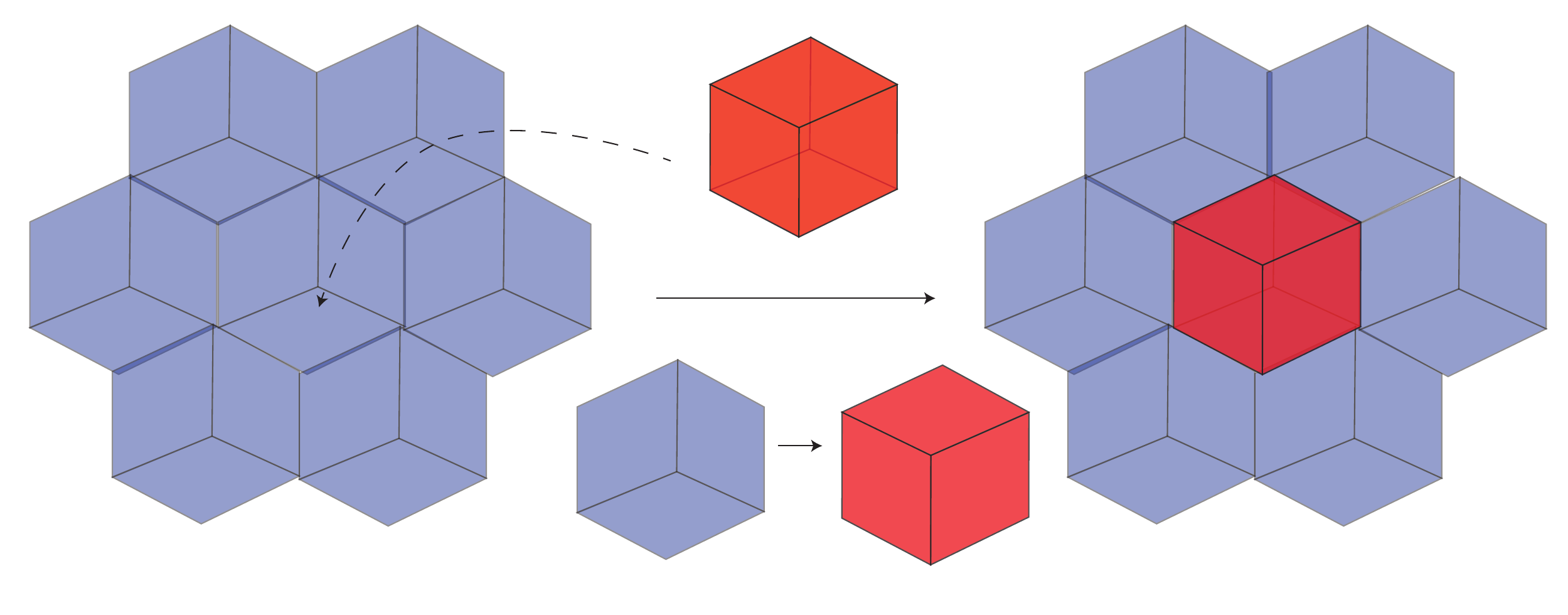}

\end{center}
 \caption{A flip and a stacking cube}\label{fig:cubes-tiling-invert}
\end{figure}

Let us suppose that cubes are stacked in the way shown in the left of Fig. \ref{fig:cubes-tiling-invert}.
If one puts a (red) cube on the top of the others, this inverts a tiling of one octagon as shown in Fig. \ref{fig:cubes-tiling-invert}.
In \cite{HenriquesSpeyer} it is shown that the flips of rhombile tilings of the planar figure satisfy two interesting relations, which we call {\em the far-commutativity} and {\em the octagon relation}. In \cite{HenriquesSpeyer} it is proved that all sequences of flips such that the initial and final tilings are same, are equivalent up to the far-commutativity and the octagon relations. 

\begin{figure}[h!]
\begin{center}
 \includegraphics[width =10cm]{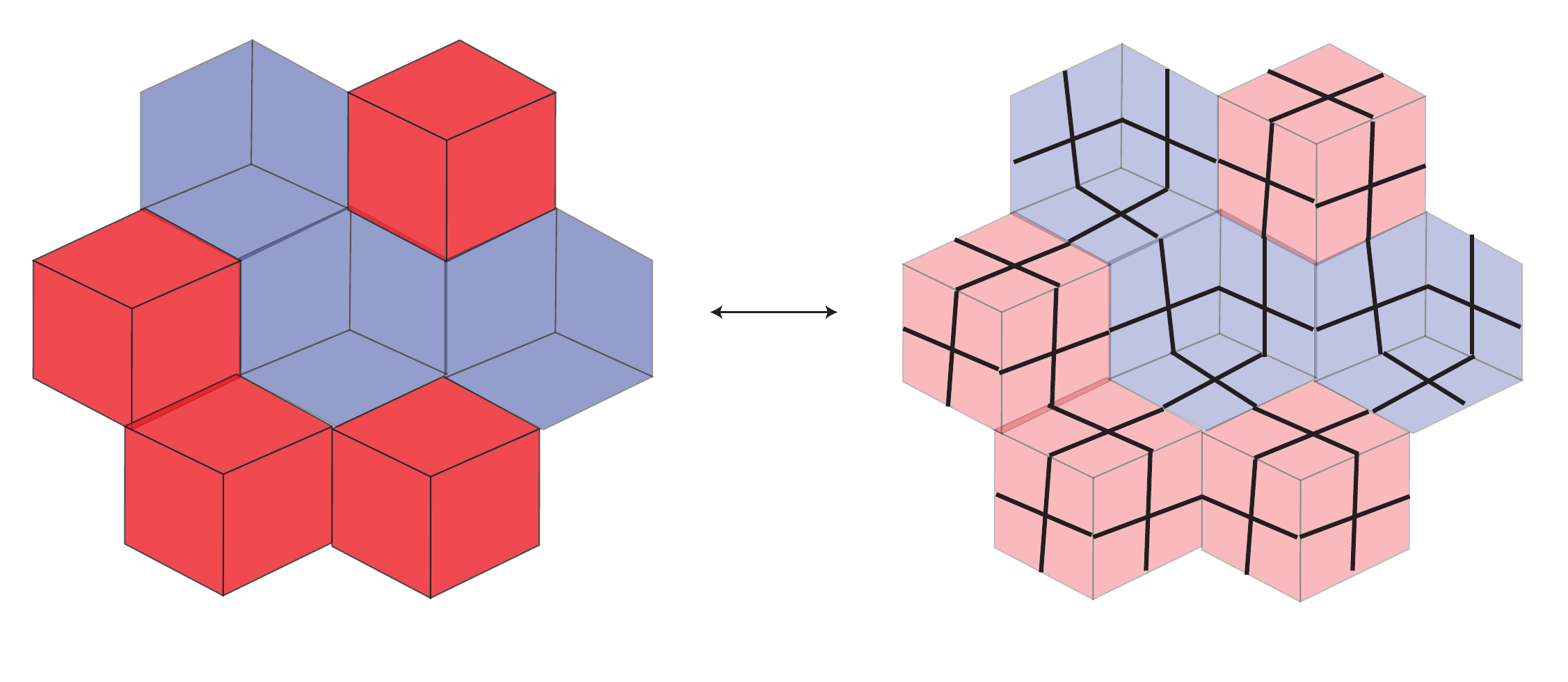}
\end{center}
 \caption{Rhombile tiling and dual diagram}\label{fig:cubes-tiling-dual}
\end{figure}

\begin{figure}[h!]
\begin{center}
 \includegraphics[width =10cm]{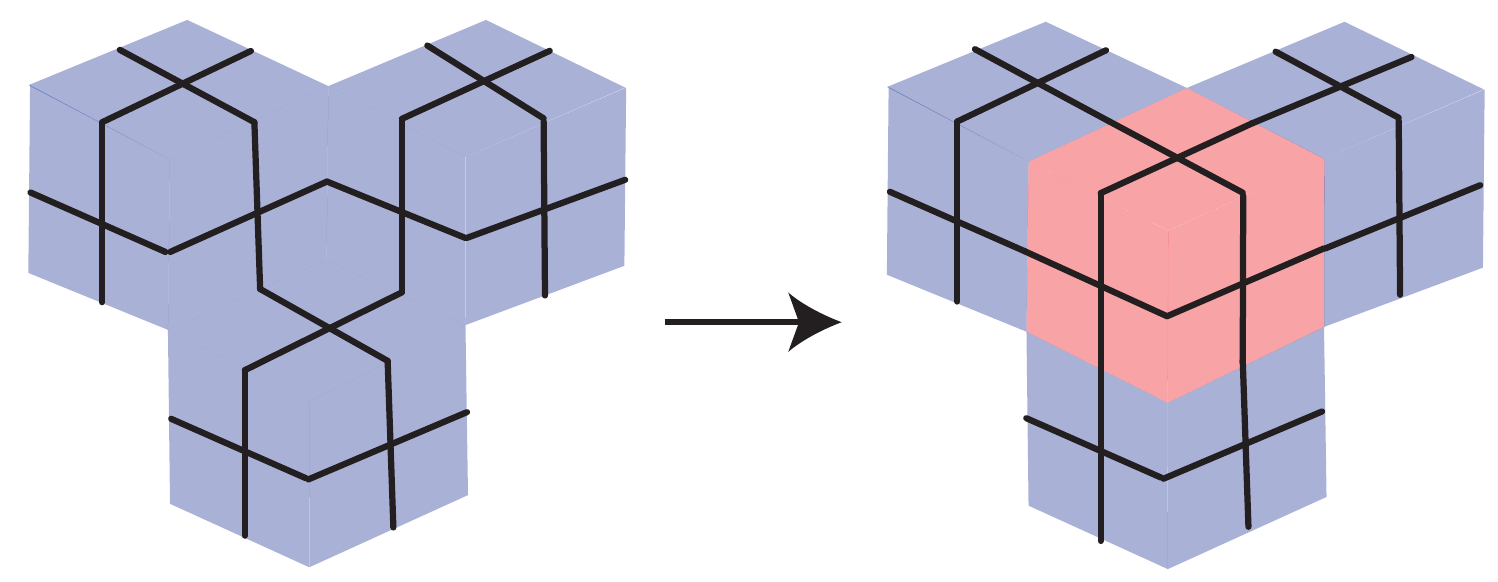}

\end{center}
 \caption{A flip and the corresponding move for dual diagrams}\label{fig:inversion-dual}
\end{figure}

If we pass to the dual diagram, this will look as shown in Fig. \ref{fig:cubes-tiling-dual}. 
In particular, the reader easily recognises that the stacking of one cube or inversion of tilings corresponds to the third Reidemeister move as shown in Fig. \ref{fig:inversion-dual}.

\begin{figure}[h!]
\begin{center}
 \includegraphics[width =6cm]{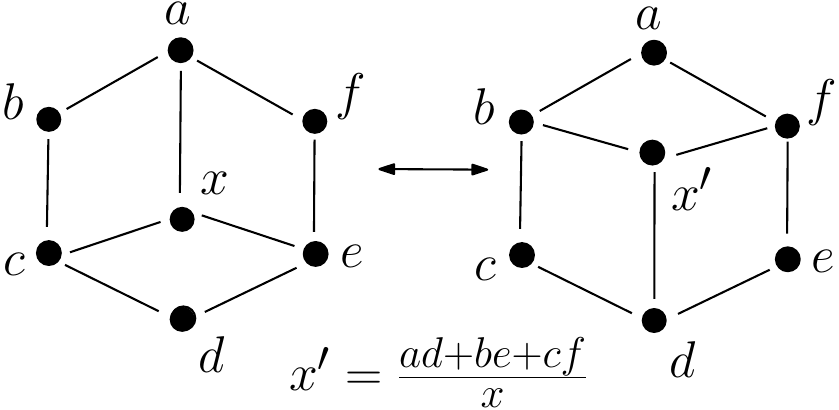}

\end{center}
 \caption{A formula associated with an inversion}\label{fig:HS-rel}
\end{figure}

 In \cite{HenriquesSpeyer} it is shown how to relate rhombile tilings to a cluster algebra. With each vertex of a tiling $T$ we associate a variable $x_{i}$, $i\in I$, see Fig.~\ref{fig:HS-rel}. 
 Assume that a tiling $T'$ is obtained from $T$ by one inversion. Say, a vertex with a variable $x$ in $T$ inside an hexagon with vertices $a,b,c,d,e,f$ is replaced by a vertex with a variable $x'$. Then $x'$ is determined by the formula
 \begin{equation}\label{HS-formula}
      x'=\frac{ad+be+cf}{x},
 \end{equation}
 and this formula provides the structure of a cluster algebra seed mutation, for details, see \cite{HenriquesSpeyer}. 

 So, following \cite{Manturov-braction}
one may expect that beautiful pictures and theorems concerning configurations of lines and
points give rise to the fact that some transformations satisfy certain equations, for example, the well known theorem about cubical curves sounds as follows.

\begin{proposition}
    Let $C, C_{1},C_{2}$ be three cubics in $\mathbb{R}P^{2}$. If $C$ goes through eight of the nine intersection points of $C_{1},C_{2}$ then $C$ goes through the ninth intersection point as well.
\end{proposition}


The reader familiar with
\cite{FominPylyavskyy} will recognise the Desargues flips on configurations of lines and points.
One of the authors (V.O.M) defined a series of groups $G_{n}^{k}$ for $n>k$ \cite{M22, Ma,MNGn3}. Those groups may be regarded as a certain generalisation of braid groups and dynamical systems led to the discovery
of the following fundamental principle:

{\em If dynamical systems describe a motion of $n$ particles and there exists a good codimension $1$ property governed by exactly $k$ particles
then these dynamical systems admit an invariant valued in the group $G_{n}^{k}$}.

In \cite{MNGn3} a partial case of this general principle was calculated explicitly: when considering a motion
of $n$ pairwise distinct points on the plane and choosing the generic codimension 1 property to be ``three
points are collinear'' or ``four points lie on a circle'', we get a homomorphism from the $n$-strand pure braid group $PB_{n}$ to the group $G_{n}^{3}$ or $G_{n}^{4}$.
In particular, the generators $a_{ijk}$ of $G_{n}^{3}$ correspond to the moments when three points are collinear. 

Here we would like to emphasise that particles need not be just points. In the present paper we consider the ``dual'' approach when we consider moving (projective) lines on the (projective) plane instead of points.

The groups $G_{n}^{k}$ are large enough. One way to show that the groups $G_{n}^{k}$ are rich, is to construct MN-indices in \cite{MNGn3}. As it was mentioned in the paper \cite{MNGn3} there are lots of homomorphisms from $G_{n}^{k}$ to free products of cyclic groups due to the so-called MN-indices. These homomorphisms allow one to extract many powerful invariants of braids and their generalisations.
Now it becomes clear that each of the $G_{n}^{k}$ theories should have its {\em dual}: instead of points in the projective space one can consider dual hyperplanes.
With a bit of phantasy, one can generalise this approach to some curvilinear objects.
The case when we pass from points to pseudo-lines is described in detail in Chapter 13.2 in \cite{FKMN}.

In the present paper it is proved that there exists a map from sequences of flips on rhombile tilings to $G_{n}^{3}$. The present paper is organised as follows:

In Section \ref{sec:rhombi-cube}, we introduce rhombile tilings and flips, and their properties. In Section \ref{sec:realisable}, we describe relationship between flips on rhombile tilings of a zonogon with $2n$ vertices and elements of group $G_{n}^{3}$. In Section \ref{sec:rhombi-RM3}, we show that flips on rhombile tilings of a zonogon with $2n$ vertices can be interpreted by lines on the zonogon and the third Reidemeiter moves. In Section \ref{sec:rhombi-nonoriented}, analogously we study the relationship between flips on rhombile tilings on the projective plane $\mathbb{R}P^{2}$ and elements of group $G_{n}^{3}$. It will be shown that there exists a sequence of rhombile tilings on $\mathbb{R}P^{2}$ such that its initial and final tilings are same, but it is not equivalent to one element sequence up to relations of flips (Example \ref{exa:nontrivial-el}). In Section \ref{sec:conf-plane-line}, we consider another interpretation of rhombile tiling and flips by means of the Desargues theorem.
In Section \ref{sec:Gnk-MN} we introduce the one of tools to study $G_{n}^{3}$ called Manturov-Nikonov indices. In Section \ref{sec:further} we formulate our further research directions.

\section{rhombile tilings, flips and stacking cubes}\label{sec:rhombi-cube}
In this section we introduce the notions of rhombile tilings and flips. We borrow the definitions from  \cite{HenriquesSpeyer}.

\begin{figure}[h!]
\begin{center}
 \includegraphics[width =8cm]{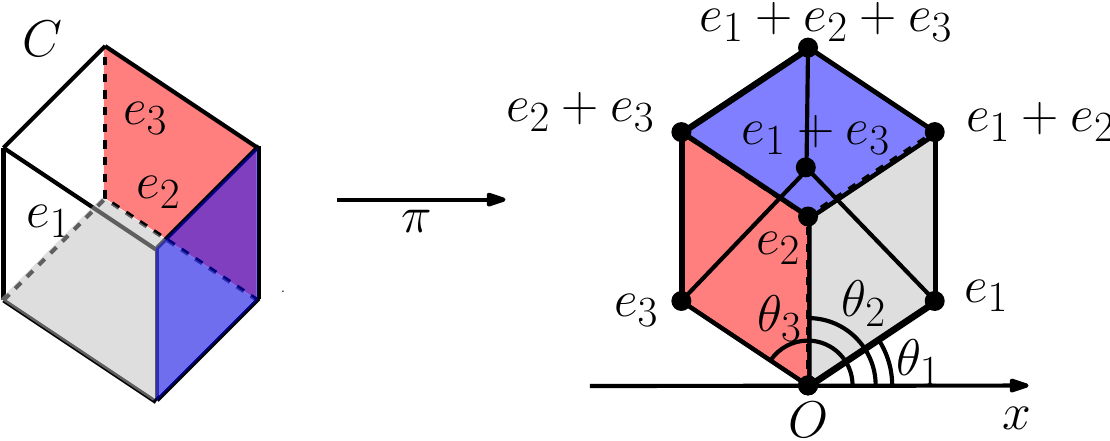}

\end{center}
 \caption{$C=\Pi_{i=1}^{3}[0,1]$ and its image $P=\pi(C)$ on $\mathbb{R}^{2}$ by $\pi: C \rightarrow \mathbb{R}^{2}$.}\label{fig:exa-projection-zonogon}
\end{figure}

Let
$$C:=\Pi_{i=1}^{n}[0,1] = \{x_{1}e_{1}+\dots+x_{n}e_{n}| 0 \leq x_{i} \leq 1\},$$
where $\{e_{1},\cdots, e_{n}\}$ is the standard basis of $\mathbb{R}^{n}$. We choose $n$ real numbers $\theta_{1},\dots, \theta_{n}$ satisfying $0<\theta_{1}<\cdots<\theta_{n}<\pi$ and let $v_{i}:=(cos\theta_{i},sin\theta_{i})\in \mathbb{R}^{2}$. Let $\pi : C \rightarrow \mathbb{R}^{2}$ be the map $x_{1}e_{1}+\dots+x_{n}e_{n} \mapsto \Sigma x_{i} v_{i}$, where $0 \leq x_{i} \leq 1$. 
Let $P=\pi(C)$. The image $P$ of $C$ is the convex hull of the images of $x_{1}e_{1}+\dots+x_{n}e_{n}$, where $x_{i} \in \{0,1\}$ for $i =1,\dots, n$. The boundary of $P$ is the polygon $Z$ having $2n$ vertices, $\pi(0)$, $\pi(e_{1})$, $ \pi(e_{1}+e_{2})$, $\dots$, $ \pi(e_{1}+e_{2}+\dots+e_{n-1})$, $\pi(e_{1}+e_{2}+\dots+e_{n-1}+e_{n})$, $ \pi(e_{2}+\dots+e_{n-1}+e_{n})$, $\dots$, $ \pi(e_{n-1}+e_{n})$, $\pi(e_{n})$ such that the $i$-th and $(n+i)$-th edges are parallel and of the same length $1$. A polygon whose edges have this property is called {\em a zonogon}.
See Fig.~\ref{fig:exa-projection-zonogon}.
{\em A tiling of $P$} is a two dimensional sub-complex $T \subset C$ such that $\pi: T\rightarrow \pi(T) \subset P$ is a homeomorphism. 

\begin{proposition}[\cite{HenriquesSpeyer}]
The map $\pi:C\rightarrow P$ induces a bijection between tilings $T \subset C$ and decompositions $\mathcal{T} =  \pi(T)$ of $Z$ into rhombi with side length $1$. 
\end{proposition}
In this paper, decompositions of $Z$ into rhombi with side length $1$ are called {\em rhombile tilings of the zonogon $Z$}. We use tilings of $C$ together with rhombile tilings of $Z$ indiscriminately.

Let $I$ be a vertex of $C$ and let $0< j<k<l\leq n$. Assume that the cube
$$c= \{ I+ xe_{j}+ye_{k}+ze_{l}~|~x,y,z\in[0,1]\}$$
is contained in $C$. The three faces of $c$ containing $I+e_{k}$ are called the {\em bottom faces of $c$} and those containing $I+e_{j}+e_{l}$ the {\em top faces}. If a tiling $T$ contains the top faces of $c$, then the complex obtained by replacing the top faces of $c$ by the bottom faces is also a tiling denoted by $T_{c}$. The converse is also possible. This operation is called {\em a flip}.

\begin{remark}
    Let $T$ be a tiling and let $c$ be the cube described above. One can interpret the flip as a stacking of a cube on the tiling $T$ as described in Fig. \ref{fig:cubes-tiling-invert}. In natural way, one can obtain a cobordism between $T$ and $T_{c}$, that is, 3-dimensional manifold with boundaries $T$ and $T_{c}$. We call $T$ and $T_{c}$ {\em the upper boundary} and {\em the lower boundary}.
\end{remark}

\begin{remark}
    Note that for a given cube $c$, we obtain two flips: bottom faces to top faces and the opposite. If $T_{c}$ is the tiling obtained from $T$ by applying the flip derived from $c$, then $(T_{c})_{c}$ and $T$ provide the same tiling.
\end{remark}

The following Lemmas can be easily proved:
\begin{lemma}[Far-commutativity]\label{lem:far-comm}
    Let $T$ be a tiling of $C$. Let $T_{c}$ be the tiling obtained from $T$ by the flip corresponding to the cube $c= \{ I+ xe_{j}+ye_{k}+ze_{l}~|~x,y,z\in[0,1]\}$. Let $c=\{ I+ xe_{j}+ye_{k}+ze_{l}~|~x,y,z\in[0,1]\}$ and $c'=\{ I+ xe_{s}+ye_{t}+ze_{u}~|~x,y,z\in[0,1]\}$ be two cubes such that $|\{e_{j},e_{k},e_{l}\}\cap \{e_{s},e_{t},e_{u}\} | \leq 1$. Then $(T_{c})_{c'}$ and $(T_{c'})_{c}$ provide the same tiling, see Fig.~\ref{fig:inversion-farcomm}.
\end{lemma}

\begin{lemma}[Octagon relation]\label{lem:octagonrel}
    Let $t = \{ I + xe_{i}+ye_{j}+ze_{k}+we_{l}~|~x,y,z,w\in[0,1]\}$ be a tesseract which includes four cubes $c_{ijk}$, $c_{ijl}$, $c_{ikl}$, $c_{jkl}$, where $c_{stu}$ is the cube generated by $e_{s}$, $e_{t}$, $ e_{u}$. Let $T$ be a tiling. Then two tilings $(((T_{c_{ijk}})_{c_{ijl}})_{c_{ikl}})_{c_{jkl}}$ and $(((T_{c_{jkl}})_{c_{ikl}})_{c_{ijl}})_{c_{ijk}}$ provide the same tiling, see Fig.~\ref{fig:octagon_rel}.
\end{lemma}

 The octagon relation in Lemma~\ref{lem:octagonrel} is derived from the 3-dimensional cube. In the next section we will show that the octagon relation is closely related to the group $G_{n}^{3}$. In general, one can obtain analogous relations from higher dimensional cubes and they must be related with the groups $G_{n}^{k}$, $k>3$.

The group $G_{n}^{k}$ and the map from the braid group to $G_{n}^{3}$ enjoy many nice properties, in particular, as shown in \cite{M22}.

\section{The realisable counterpart of $G_{n}^{3}$}\label{sec:realisable}

Assume that we are stacking cubes starting from tilings,  which are 2-dimensional pictures, and look at the picture ``on the top''. 
After two different sequences of stacking of cubes we can get the same tiling as asserted in Lemma \ref{lem:far-comm} and \ref{lem:octagonrel}. This fact is closely related to the structure of groups $G_{n}^{3}$ and their relations. Different ways of applying generators of $G_{n}^{3}$ may lead to the same result.

In our situation, a cube being stacked corresponds to a generator of the group $G_{n}^{3}$, which one can associate with indices, so-called {\em MN-indices}. One can expect that MN-indices reflect the summation of volumes of stacked cubes obtained from a sequence of flips.

 \begin{definition}
Let $n,k$ be integers such that $1<k\leq n$. The group $G_{n}^{k}$ is given by the group presentation with generators $\{ a_{m} ~|~ m \subset \{1, \cdots, n \}, |m| = k \}$ and relations as follows:

\begin{enumerate}
\item $(a_{m})^{2} = 1$, for $ m \subset \{1, \cdots, n\}$,
\item $a_{m}a_{m'} = a_{m'}a_{m}$ for $m,m' \subset \{1, \cdots, n\}$, $|m \cap m'| <k-1$,
\item $a_{m^{1}} \cdots a_{m^{k+1}} = a_{m^{k+1}} \cdots a_{m^{1}}$, where $U = \{i_{1}, \cdots, i_{k+1} \} \subset \{1, \cdots, n \} $ and $m^{l} = U \backslash \{i_{l}\}$.
\end{enumerate}
\end{definition}

As a partial case the group $G_{n}^{3}$ can be defined as follows:
\begin{definition}
The group $G_{n}^{3}$ is given by the group presentation with generators $\{ a_{\{i,j,k\}}~|~ \{i,j,k\} \subset \{1, \dots, n\}, |\{i,j,k\}| = 3\}$ and relations as follows:

\begin{enumerate}
\item $a_{\{ijk\}}^{2} = 1$ for all $\{i,j,k\} \subset \{1, \cdots,n\}$, 
\item $a_{\{ijk\}}a_{\{stu\}} = a_{\{stu\}}a_{\{ijk\}}$, if $| \{i,j,k\} \cap \{s,t,u\} | < 2$.
\item $a_{\{ijk\}}a_{\{ijl\}}a_{\{ikl\}}a_{\{jkl\}} = a_{\{jkl\}}a_{\{ikl\}}a_{\{ijl\}}a_{\{ijk\}}$ for distinct $i,j,k,l$.
\end{enumerate}
Simply, we denote $a_{ijk}:=a_{\{ijk\}}$.
\end{definition}

\begin{figure}[h!]
\begin{center}
 \includegraphics[width =6cm]{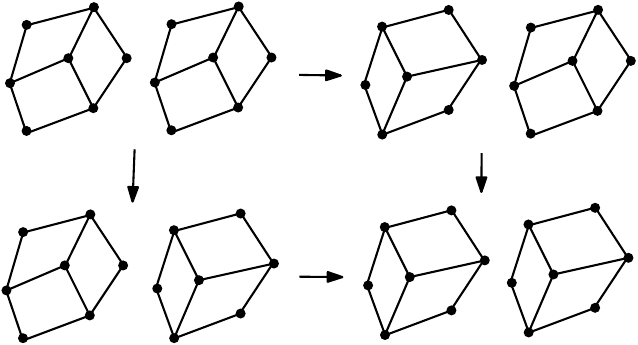}

\end{center}
 \caption{Far-commutativity of flips}\label{fig:inversion-farcomm}
\end{figure}

\begin{figure}[h!]
\begin{center}
 \includegraphics[width =6cm]{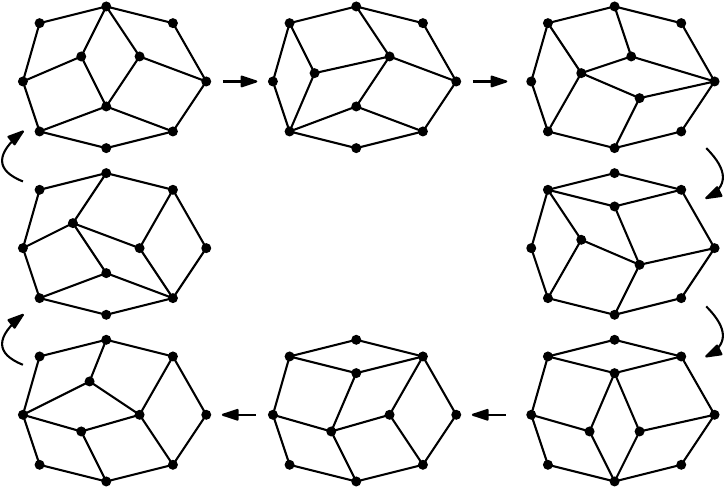}

\end{center}
 \caption{Octagon relation of flips}\label{fig:octagon_rel}
\end{figure}

Let $Z$ be a zonogon with $2n$ vertices. Let $\mathbb{X}_{1}$ be the graph whose vertices correspond to the rhombile tilings of $Z$, and whose edges are the flips. That is, \begin{eqnarray*}
   V(\mathbb{X}_{1}) &=& \{\text{rhombile tilings}~T~ \text{of}~Z\},\\
    E(\mathbb{X}_{1}) &=& \{TT_{c}\},
\end{eqnarray*}
where $T_{c}$ is the rhombile tiling obtained from $T$ by one flip on some cube $c$. One can compare $\mathbb{X}_{1}$ with the notion of a {\em flip graph} \cite{flipgraph}.
\begin{proposition}[\cite{HenriquesSpeyer}]\label{prop:conn}
        The graph $\mathbb{X}_{1}$ is connected.
\end{proposition}
We now construct a  2-dimensional cell complex $\mathbb{X}_{2}$ by gluing 2-cells into $\mathbb{X}_{1}$ as follows:
\begin{enumerate}
    \item If $T \sim T_{c}$ and $T \sim T_{c'}$ are flips involving disjoint sets of rhombi, then these two flips can be applied simultaneously to get a fourth tiling $T_{cc'}=T_{c'c}$ by far-commutativity. The vertices $T,T_{c},T_{c'},T_{cc'}$ form a 4-cycle in $\mathbb{X}_{1}$ on which we attach a square, see Fig.~\ref{fig:inversion-farcomm}.
    \item  Suppose that we take a rhombile tiling $T$ contained in one of the pictures shown in Fig. \ref{fig:octagon_rel}. Then one can perform the corresponding cycle shift of eight flips on $T$. In such a case, glue an octagon via its boundary to this series of eight flips.

\end{enumerate}

\begin{proposition}[\cite{HenriquesSpeyer}]\label{prop:simply-conn}
    The cell complex $\mathbb{X}_{2}$ is simply connected.
\end{proposition}

For an arbitrary path $T \rightarrow T'$ in $\mathbb{X}_{1}$ by Proposition \ref{prop:conn} one can find the sequence of flips on cubes $c_{1}, c_{2},\dots, c_{s}$.
Let us define an equivalence relation $\sim$ on $Path(\mathbb{X}_{1})$: Two paths are {\it equivalent} if they are homotopic relative to end points in $\mathbb{X}_{2}$. Let us define a map $\phi$ from the set $Path(\mathbb{X}_{1})/\langle\sim\rangle$ of equivalence classes of paths to $G_{n}^{3}$ as follows:
For a path $T \rightarrow T'$ let us define $$\phi(T\rightarrow T')= a^{c_{1}}\dots a^{c_{n}},$$
where $a^{c_{s}}=a_{ijk}$ for $c_{s}= \{ I+ xe_{i}+ye_{j}+ze_{k}~|~x,y,z\in[0,1]\}$.

\begin{theorem}
    The map $\phi : Path(\mathbb{X}_{1})/\langle\sim\rangle \rightarrow G_{n}^{3}$ is well-defined.
\end{theorem}

\begin{proof}
    Let $T$ and $T'$ be two rhombile tilings of $Z$. By Proposition \ref{prop:simply-conn} two paths from $T$ to $T'$ are homotopic in $\mathbb{X}_{2}$. Since a 2-cell of $\mathbb{X}_{2}$ is a quadrilateral or an octagon, two paths from $T$ to $T'$ bound quadrilaterlas and octagons. By Lemma \ref{lem:far-comm} and Lemma \ref{lem:octagonrel}, one can show that if two paths from $T$ to $T'$ bound one quadrilateral or one octagon, then the images of them by $\phi$ are equivalent in $G_{n}^{3}$.
\end{proof}

Let us call elements $ \beta \in \phi(Path(\mathbb{X}_{1})/\langle\sim\rangle) \subset G_{n}^{3}$ {\em realisable elements} and $\phi(Path(\mathbb{X}_{1})/\langle\sim\rangle)$ a {\em realisable counterpart of $G_{n}^{3}$}\label{def:realisable counterpert}. It is clear that $1 \in \phi(Path(\mathbb{X}_{1})/\langle\sim\rangle) \subset G_{n}^{3}$ and for any realisable element $\beta \in \phi(Path(\mathbb{X}_{1})/\langle\sim\rangle)$, there exists $\beta^{-1} \in \phi(Path(\mathbb{X}_{1})/\langle\sim\rangle)$. For two paths $T\rightarrow T'$ and $T'\rightarrow T''$, $\phi(T\rightarrow T' \rightarrow T'') = \phi(T\rightarrow T')\phi(T\rightarrow T')$. But for two elements in $\phi(Path(\mathbb{X}_{1})/\langle\sim\rangle)$ one cannot guarantee that their corresponding paths in $Path(\mathbb{X}_{1})$ can be composed. In the same time it is not known that the equivalence relation in $\phi(Path(\mathbb{X}_{1})/\langle\sim\rangle)$ is induced from the equivalence relation in $G_{n}^{3}$. Therefore we cannot guarantee that $\phi(Path(\mathbb{X}_{1})/\langle\sim\rangle)$ is a subgroup, but we consider it as a ``subset'' of $G_{n}^{3}$. 

If we consider the subset $RG_{n}^{3} \subset \phi(Path(\mathbb{X}_{1})/\langle\sim\rangle)$ of images of all closed paths containing a fixed vertex $T_{0}$ in $\mathbb{X}_{1}$, then $RG_{n}^{3}$ is a subgroup of $G_{n}^{3}$. Let us call $RG_{n}^{3}$ {\em a realisable subgroup of $G_{n}^{3}$ by rhombile tilings on a zonogon $Z$ with $2n$ vertices}.
\begin{corollary}\label{cor:trivial-from-plane}
    For a closed path $\beta \in Path(\mathbb{X}_{1})/\langle\sim\rangle$, $\phi(\beta) = 1$ in $G_{n}^{3}$
\end{corollary}
\begin{proof}
    This corollary follows from Proposition~\ref{prop:simply-conn}.
\end{proof}

That is, the realisable subgroup of $G_{n}^{3}$ by rhombile tilings and flips on $Z$ is trivial.

\section{Rhombile tilings and the third Reidemeister moves}\label{sec:rhombi-RM3}

Interestingly, flips of rhombile tilings are related to Reidemeister moves in a natural way:
For each rhombus of a rhombile tiling $T$ corresponding to $\{I+xe_{i}+ye_{j}~|~x,y \in [0,1]\}$ let us take arcs $I+xe_{i}+\frac{1}{2}e_{j}$ and $I+\frac{1}{2}e_{i}+ye_{j}$ for $x\in [0,1]$ 
 and $y\in [0,1]$. They intersect at $I+\frac{1}{2}e_{i}+\frac{1}{2}e_{j}$. The set of arcs obtained from a rhombile tiling as above is called {\em the dual diagram}.
 
\begin{figure}[h!]
\begin{center}
 \includegraphics[width =5cm]{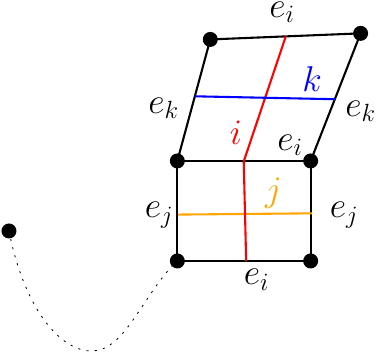}

\end{center}
 \caption{Dual diagram of a rhombile tiling and labels for arcs}\label{fig:dual-rhombi}
\end{figure}

 For each arc, we associate the label $i$ to the arc $I+\frac{1}{2}e_{i}+ye_{j}$, $y\in [0,1]$. By definition of a dual diagram, the label is defined compatibly, for example, see Fig.~\ref{fig:dual-rhombi}. 

\begin{figure}[h!]
\begin{center}
 \includegraphics[width =10cm]{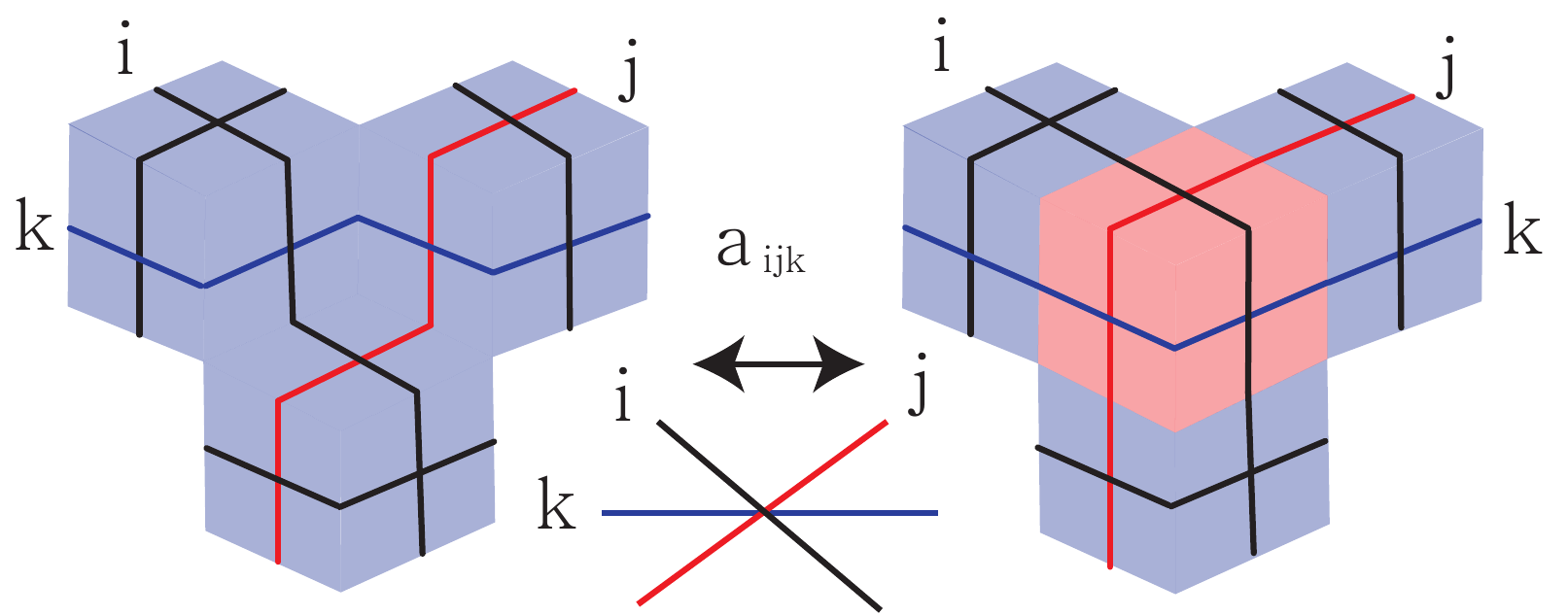}

\end{center}
 \caption{A flip, the corresponding move for dual diagrams and the associated generator of the group $G_{n}^{k}$}\label{fig:inversion-dual-Gn3}
\end{figure}

For two rhombile tilings $T$ and $T_{c}$, where $T_{c}$ is the rhombile tiling obtained from $T$ by the flip on a cube $c$, two dual diagrams corresponding to $T$ and $T_{c}$ are obtained from each others by the third Reidemeister move, see Fig.~\ref{fig:inversion-dual-Gn3}.

\begin{figure}[h!]
\begin{center}
 \includegraphics[width =10cm]{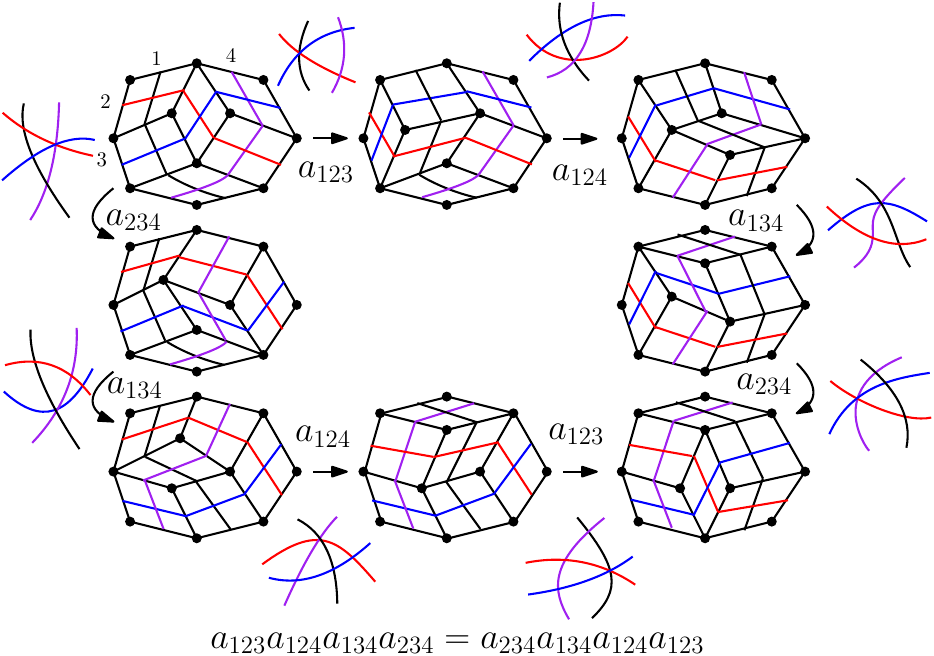}

\end{center}
 \caption{Octagon relation of dual diagrams from tilings and the relation of $G_{n}^{3}$}\label{fig:octagon_rel-dual-Gn3}
\end{figure}

Note that applying the third Reidemeister move a triple point of three arcs appears.  
Now we associate a triple point of arcs numbered by $i,j,k$ to a generator $a_{ijk}\in G_{n}^{3}$. From Fig.~\ref{fig:inversion-farcomm} and \ref{fig:inversion-dual-Gn3} it is easy to see that sequences of the third Reidemeister moves provide relations $a_{ijk}^{2}=1$ and $a_{ijk}a_{stu}= a_{stu}a_{ijk}$ for $|\{i,j,k\} \cap \{s,t,u\}|<2$. From the octagon relation for the third Reidemeister moves we obtain $a_{ijk}a_{ijl}a_{ikl}a_{jkl} = a_{jkl}a_{ikl}a_{ijl}a_{ijk}$ for distinct $i,j,k,l$ as described in Fig.~\ref{fig:octagon_rel-dual-Gn3}.

\section{Rhombile tiling of $\mathbb{R}P^{2}$ and realisable elements}\label{sec:rhombi-nonoriented}

In \cite{Manturov-braction} the second author studied actions of braids on the configuration space of the points in $\mathbb{R}P^{2}$. In the present section, we define the rhombile tilings on the projective plane $\mathbb{R}P^{2}$ by using rhombile tilings on a zonogon $Z$ with $2n$ vertices and construct a map from any sequence of flips to $G_{n}^{3}$. One might expect the construction to be similar to that of the plane. However, unlike rhombile tilings on the plane, it gives rise to nontrivial realisable elements in $G_{n}^{3}$:

Rhombile tilings on $\mathbb{R}P^{2}$ and the Klein bottle are obtained as follows:
let redraw the 2n-gon as a rectangle. Let us assume that the edges on the boundary are labeled by numbers $e =(e_{1},\dots, e_{n})$ where $e_{i} \in \{1,\dots,n\}$ such that $e_{i} \neq e_{j}$ for $i\neq j$. Note that for a given rhombile tiling of the 2n-gon the labels of edges of the 2n-gon can be induced from the labels $e$ of boundary edges as the case of a zonogon with $2n$ vertices.

\begin{figure}[h!]
\begin{center}
 \includegraphics[width =5cm]{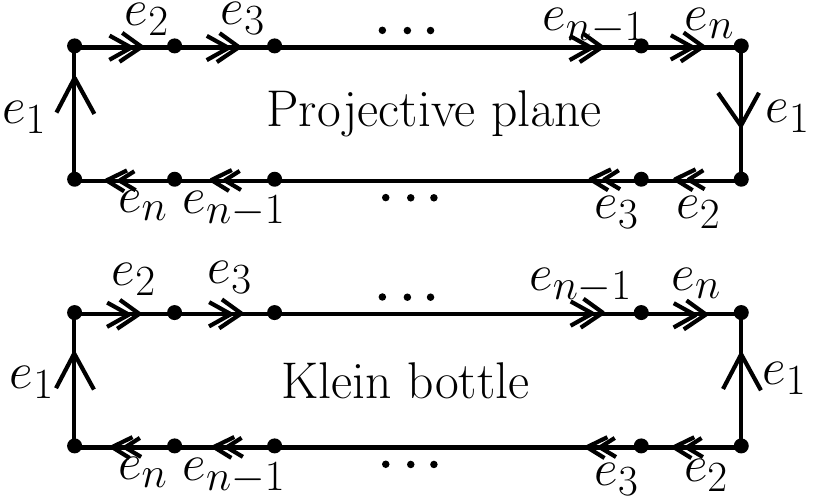}

\end{center}
 \caption{The projective plane and the Klein bottle}\label{fig:2ngon-nonori}
\end{figure}

By attaching boundaries as described in Fig.~\ref{fig:2ngon-nonori} we obtain $\mathbb{R}P^{2}$ and the Klein bottle. 
Let $S^{e}$ be the projective plane $\mathbb{R}P^{2}$ or the Klein bottle obtained from the 2n-gon with labels $e =(e_{1},\dots, e_{n})$ on boundary edges. 

Since we can glue the surface $S^{e}$ by identifying some pairs of sides of the 2n-gon with labels $e$, we may obtain rhombile tilings of $S^{e}$ from rhombile tilings of the 2n-gon with labels $e$. In particular, for a given rhombile tiling $T^{e}$ of $S^{e}$ the labels of the edges of $T^{e}$ are induced from the labels $e$. From now on, we assume that tilings $T^{e}$ of $S^{e}$ have labels of edges induced from $e$.

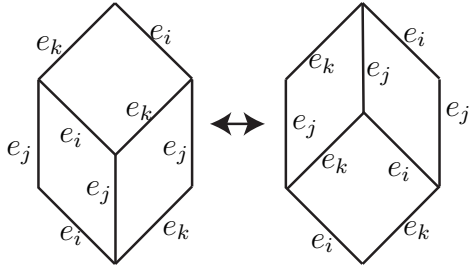
\begin{figure}[h!]
\begin{center}
 \includegraphics[width =6cm]{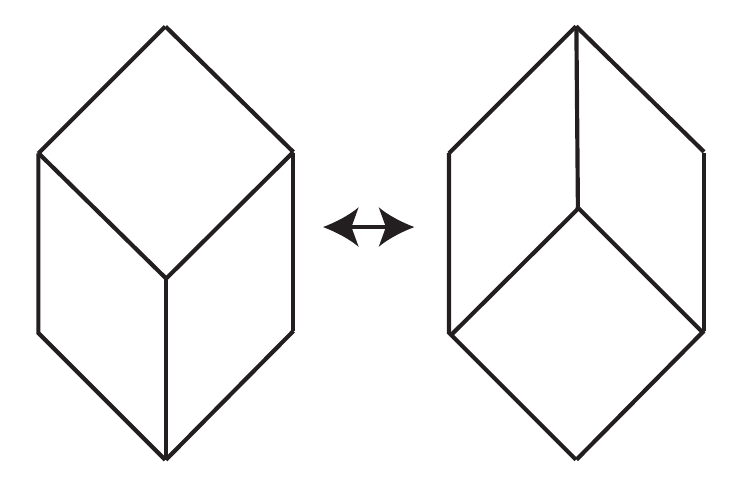}
\put(-154, 18){$e_{i}$}
\put(-154, 55){$e_{i}$}
\put(-120, 94){$e_{i}$}
\put(-144, 35){$e_{j}$}
\put(-174, 50){$e_{j}$}
\put(-115, 50){$e_{j}$}
\put(-115, 20){$e_{k}$}
\put(-129, 65){$e_{k}$}
\put(-163, 90){$e_{k}$}
\put(-59, 15){$e_{i}$}
\put(-24, 93){$e_{i}$}
\put(-30, 42){$e_{i}$}
\put(-66, 60){$e_{j}$}
\put(-8, 65){$e_{j}$}
\put(-38, 80){$e_{j}$}
\put(-24, 20){$e_{k}$}
\put(-55, 45){$e_{k}$}
\put(-60, 83){$e_{k}$}
\end{center}
 \caption{A flip on an octagon $c$}\label{fig:flip-normal}
\end{figure}

For a rhombile tiling $T^{e}$, let $c$ be an hexagon consisting of three quadrilaterals. There are two possible rhombile tilings restricted in the octagon $c$ as described in Fig.~\ref{fig:flip-normal}.
{\it A flip on an hexagon $c$ of $T^{e}$} is defined by the replacement of a given rhombile tiling inside $c$ by another. Note that $c$ consists of edges labeled by three different numbers $e_{i},e_{j},e_{k}$ induced from $e$. Let us denote the tiling obtained from $T^{e}$ by the flip on $c$ by $T_{c}^{e}$.

Let $\mathbb{X}_{1}^{S}$ be a graph such that the tilings $T^{e}$ of $S^{e}$ for any label $e$ are associated to vertices and edges $T^{e}T_{c}^{e}$ are connected, when $T_{c}^{e}$ is the tiling obtained from $T^{e}$ by the flip on $c$.

\begin{corollary}
    Flips on tilings of $S$ satisfy the following relations:
    \begin{itemize}
        \item $(T_{c})_{c}=T$,
        \item Far-commutativity,
        \item Octagon relation.
    \end{itemize}
\end{corollary}

Let us fix a rhombile tiling $T$ with a fixed label $e$. Let $Path(\mathbb{X}_{1}^{S}, T)$ be the set of all paths in the graph $\mathbb{X}_{1}^{S}$ starting from $T$. Let us define an equivalence relation $\sim'$ on $Path(\mathbb{X}_{1}^{S}, T)$ as follows: Two paths are equivalent if one can obtained from another by far-commutativity and octagon relations.

Now let us define $\phi^{S} : Path(\mathbb{X}_{1}^{S},T)/\langle \sim' \rangle \rightarrow G_{n}^{3}$ by
for a path $T \rightarrow (\dots((T_{c_{1}})_{c_{2}})\dots)_{c_{n}}$ in $\mathbb{X}_{1}^{S}$,  $$\phi^{S}(T \rightarrow (\dots((T_{c_{1}})_{c_{2}})\dots)_{c_{n}}) = a^{c_{1}}\dots a^{c_{n}},$$
where $a^{c_{s}}=a_{ijk}$ for a cube $c_{s}$ consisting of $e_{i},e_{j},e_{k}$.

\begin{corollary}
    The map $\phi^{S}:Path(\mathbb{X}_{1}^{S}, T)/\langle \sim' \rangle \rightarrow G_{n}^{3}$ is well-defined.
\end{corollary}

Let us call elements $ \beta \in \phi(Path(\mathbb{X}_{1}^{S}, T)/\langle \sim' \rangle) \subset G_{n}^{3}$ {\em realisable elements by tilings of $S$} and $\phi(Path(\mathbb{X}_{1}^{S}, T)/\langle \sim' \rangle)$ {\em realisable counterpart of $G_{n}^{3}$ by rhombile tilings of $S$}. 

Similarly to the case of $\phi(Path(\mathbb{X}_{1})/\langle\sim\rangle)$ defined in page~\pageref{def:realisable counterpert} we cannot guarantee that $\phi^{S}(Path(\mathbb{X}_{1}^{S}, T)/\langle \sim' \rangle)$ is a subgroup, but the subset $R^{S}G_{n}^{3} \subset \phi(Path(\mathbb{X}_{1}^{S}, T)/\langle \sim' \rangle)$ of the image of all closed paths containing a fixed vertex $T$ is a subgroup of $G_{n}^{3}$. We call $R^{S}G_{n}^{3}$ {\em a realisable subgroup of $G_{n}^{3}$ by rhombile tilings on $S$}.

\begin{example}\label{exa:nontrivial-el}
    For $S = \mathbb{R}P^{2}$ there is a closed path $p$ in $\phi^{S}(Path(\mathbb{X}_{1}^{S}, T)/\langle \sim' \rangle)$ such that $\phi^{S}(p) \in G_{n}^{3}$ is non-trivial. For example, a closed path $p$ in $Path(\mathbb{X}_{1}^{S}, T)/\langle \sim' \rangle$ is given as described in Fig.~\ref{fig:nontrivial-real-el}. Then, by definition of $\phi^{S}$, $\phi^{S}(p) = a_{124}a_{123}a_{124}a_{123}$ in $G_{4}^{3}$ and one can show that it is nontrivial element (see Example~\ref{exa:nontriviality}).
    \begin{figure}[h!]
\begin{center}
 \includegraphics[width =8cm]{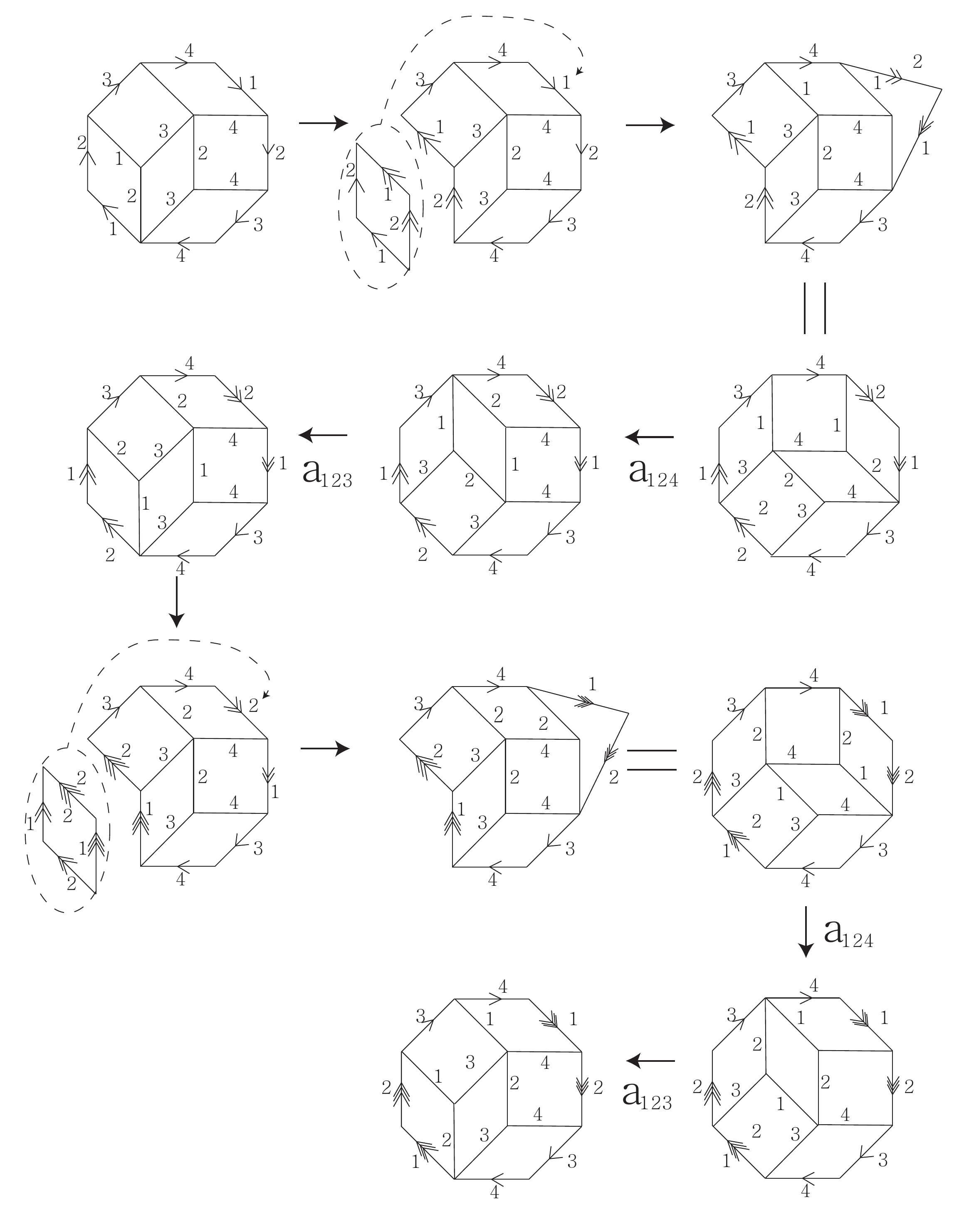}

\end{center}
 \caption{A closed path $p$ in $\phi^{S}(Path(\mathbb{X}_{1}^{S}, T)/\langle \sim' \rangle)$ such that $\phi^{S}(p) \in G_{n}^{3}$ is non-trivial.}\label{fig:nontrivial-real-el}
\end{figure}

\end{example}
Here we would like to emphasise that in Corollary~\ref{cor:trivial-from-plane} it is shown that closed paths on configuration space of lines on the plane provide the trivial element in $G_{n}^{3}$, but configurations of lines on the projective plane have nontrivial elements as images of the function.

\begin{remark}
The present construction cannot be applied to the case of orientable surfaces analogously, because when we obtain oritentable surfaces by attaching boundaries of the 4n-gon the labels of edges of rhombile tilings of $S$ cannot be derived from the labels of edges of rhombile tilings of the 4n-gon.

\end{remark}

\section{Configurations of planes and points}\label{sec:conf-plane-line}

In this section\footnote{The present section is not quite following the story of the paper. Nevertheless, we include it in order to emphasise the following. In mathematics, one and the same structure appears under different names in different areas, and this often turns out to be fruitful for all parts of mathematics where it appears. We shall touch on the relation sketched here in a separate publication. Some pictures are kindly borrowed from the paper
\cite{FominPylyavskyy}.}, we consider another way to associate configurations of lines on the projective plane to rhombile tilings based on the Desargues theorem. Our new heroe is the notions of {\bf tile} and {\bf coherent tile}. Here we just cite Fomin and Pylyavskyy \cite{FominPylyavskyy} verbatim and cite their paper in quotes.

\color{red} Let $\R{}P^{2}$ be a real projective space.
We denote by $\R{}P^{2*}$ the set of hyperplanes in~$\R{}P^{2}$.
In particular, when $\R{}P^{2}$ is a plane, the elements of $\R{}P^{2*}$ are lines.
A~point $A\in \R{}P^{2}$ and a hyperplane $\ell\in\R{}P^{2*}$ are called \emph{incident} to each other if $A\in\ell$.

We denote by $(AB)$ the line passing through two distinct points $A$ and~$B$.

A  \emph{tile} is a topological quadrilateral
(that is, a closed oriented disc with four marked points on its boundary)
whose vertices are clockwise labeled $A, \ell, B, m$,
where $A,B\in\R{}P^{2}$ are points and $\ell,m\in\R{}P^{2*}$ are hyperplanes:
\begin{equation}
\label{eq:AlBm}
\centering\includegraphics[width=60pt]{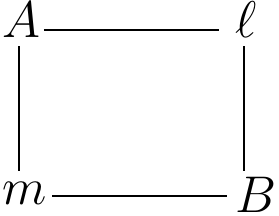}
\end{equation}
Such a tile is called \emph{coherent} if
\begin{itemize}
\item neither $A$ nor $B$ is incident to either $\ell$ or~$m$;
\item either $A=B$ or $\ell=m$ or
else the intersection of the line $(AB)$ and the codimension~$2$ subspace $\ell\cap m$ is nonempty.
\end{itemize}

\begin{figure}
\centering\includegraphics[width=120pt]{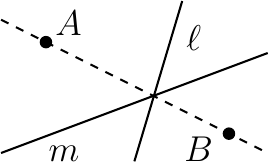}
\caption{Definition of a coherent tile} \label{fig:coherent-tile}
 \end{figure}

In the case of the projective plane, a coherent tile involves
two points $A$, $B$ and two lines $\ell$, $m$ not incident to them such that either $A=B$ or $\ell =m$ or else the line $(AB)$ passes through the point $\ell \cap m$. See Fig.~\ref{fig:coherent-tile}.

\begin{figure}
\centering\includegraphics[width=150pt]{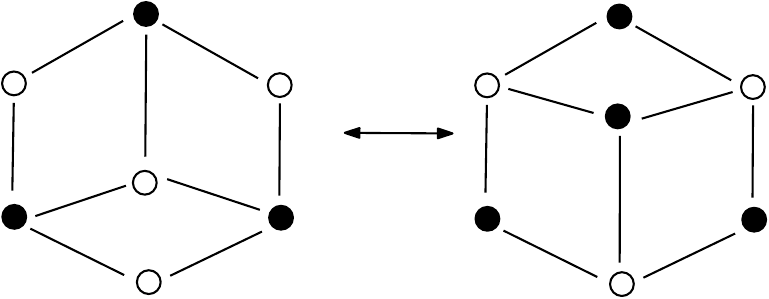}
\caption{The Desargues flip. Black dots correspond to elements of $\R{}P^{2}$, and white dots correspond to elements of $\R{}P^{2*}$.}
\label{Desarguesflip}
\end{figure}

From the Desargues Theorem, it follows that moving lines associated with a sequence of tiles such that neighbouring tiles are related by the Desargues flips described in 
 Fig. \ref{Desarguesflip}. 
Here we state the Desargues theorem.
\begin{theorem}[G. Desargues, ca. 1639]
    Let $a,b,c$ be distinct concurrent lines on the complex/real projective plane. Pick generic points $A_{1},A_{2} \in a$, $B_{1},B_{2}\in b$, $C_{1},C_{2} \in c$. Then the points $A=(B_{1}C_{1}) \cap (B_{2}C_{2})$, $B=(A_{1}C_{1})\cap (A_{2}C_{2})$, $(A_{1}B_{1})\cap(A_{2}C_{2})$ are collinear.
    \begin{figure}
\centering\includegraphics[width=150pt]{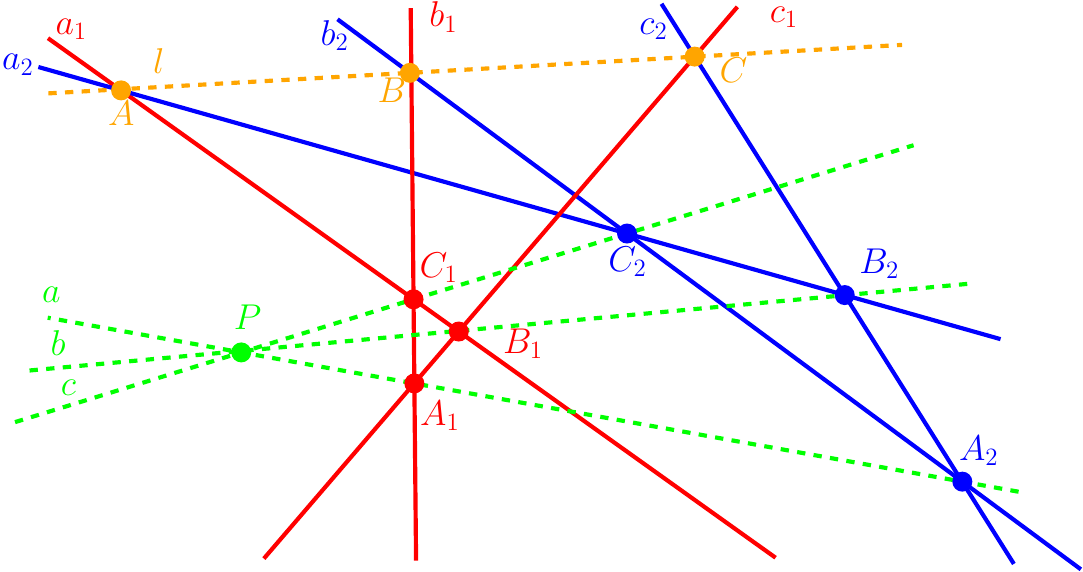}
\caption{The Desargues configuration of lines}
\label{Desargues}
\end{figure}
\end{theorem}

\begin{figure}
\centering\includegraphics[width=150pt]{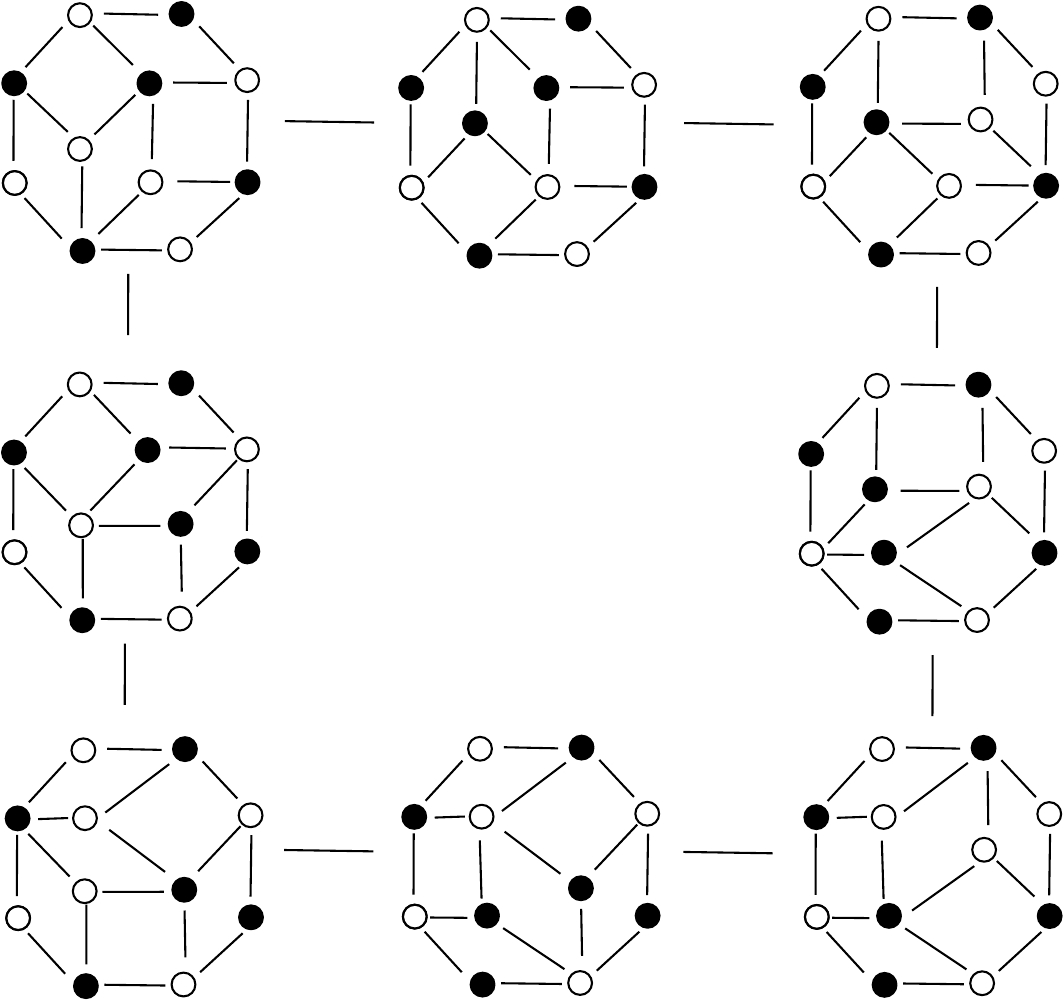}
\caption{The octagon relation}
\label{octogon1}
\end{figure}

By taking the statement from the paper by Fomin and Pylavskyy we obtain that 
the Desargues flips satisfy the octagon relation (see Fig.\ref{octogon1}).

Now one can prove that following statement:
\begin{theorem}[\cite{Manturov-braction}]
There is an action of braids on configurations of lines and planes in $\R{}P^{2}$ such that isotopic braids give rise to equal actions.
\end{theorem}

\color{black}
\section{Manturov-Nikonov indices and the group $G_{n}^{3}$}\label{sec:Gnk-MN}

The groups $G_{n}^{k}$ are interesting for a variety of reasons. One of them is that one can extract some powerful information in an easy way. Namely, one associate some ``indices'' to letters in a way such that one gets a well-defined map from the groups $G_{n}^{k}$ to free products of the group $\mathbb{Z}_{2}$. 
This makes non-triviality of many elements of $G_{n}^{k}$ visible. Here we provide the formula for the case $k=3$.

Let $\beta = F a_{ijk} G \in G_{n}^{3}$ for some $F,G \in G_{n}^{3}$. Let us denote the generator $a_{ijk}$ in $\beta$ by $c$. For $ l \in \{1,2,\cdots, n\} \backslash \{i,j,k\}$, define $i_{c}(l)$ by 
$$i_{c}(l) = (N_{jkl}+N_{ijl},N_{ikl}+N_{ijl}) \in \mathbb{Z}_{2} \times \mathbb{Z}_{2},$$
where $N_{ikl}$ is the number of the occurrence of $a_{ikl}$ from the start of $\beta$ to the letter $c$. 

\begin{example}\label{exa:nontriviality}
Let us consider an element $\beta = a_{123}a_{134}a_{145}a_{123}\underline{a_{124}}a_{134}a_{124} \in G_{5}^{3}$. Let us consider $c=a_{124}$. Then for $l=3,5$ one obtains
\begin{eqnarray*}
i_{c}(3) = (N_{134}+N_{123}, N_{234}+N_{123}) = (1,0)~\text{mod}~2,\\
i_{c}(5) = (N_{145}+N_{125}, N_{345}+N_{125}) = (1,0)~\text{mod}~2.
\end{eqnarray*}
Denote $a_{124}$ with indices by $a_{124}^{((1,0)_{3}(1,0)_{5})}$.


\end{example}

Let us fix $i,j,k\in \{1,\dots, n\}$. Let
 $$F_{n}^{3} = \langle\{ \sigma ~|~ \sigma : \{1,2,\cdots n\} \backslash \{i,j,k\} \rightarrow \mathbb{Z}_{2} \times \mathbb{Z}_{2} \}|\{ \sigma^{2} = 1\} \rangle \cong (\mathbb{Z}_{2})^{*2^{2(n-3)}}.$$ Then $i_{c}$ can be considered as a generator of $F_{n}^{3}$. 
For an element $\beta \in G_{n}^{3}$ such that $\beta$ contains $m$ generators $a_{ijk}$. Let us denote $a_{ijk}$ contained in $\beta$ as below:
\begin{eqnarray*}
    \beta &=& \cdots a_{ijk} \cdots a_{ijk} \cdots a_{ijk} \cdots a_{ijk} \cdots\\
    &=& \cdots~ c_{1}~\cdots ~c_{k}~ \cdots ~ c_{l}\cdots ~c_{m}~\cdots,
\end{eqnarray*}
where $1<k<l<m$. Define a map $w_{(i,j,k)} : G_{n}^{3} \rightarrow F_{n}^{3}$ by $$w_{(i,j,k)}(\beta) = i_{c_{1}}i_{c_{2}} \cdots i_{c_{m}}.$$

 \begin{example}

For $\beta = a_{123}a_{134}a_{123}a_{124}a_{134}a_{124}$, $w_{124}(\beta) =  i_{c_{1}}i_{c_{2}} = ((1,0)_{3}(1,0)_{5})_{c_{1}}((0,0)_{3}(1,0)_{5})_{c_{2}}$.

\end{example}

\begin{proposition}[\cite{MNGn3}]
For a positive integer $n$ and for $i, j,k \in \{1, \cdots ,n\}$ such that $|\{i,j,k \} | =3$, $w_{(i,j,k)}$ is well-defined. 
\end{proposition}

\begin{example}
    For $\beta=a_{124}a_{123}(=c_{1})a_{124}a_{123}(=c_{2})$ in $G_{4}^{3}$, we obtain
    \begin{eqnarray*}
    i_{c_{1}}(4) &=& (N_{234}+N_{124},N_{134}+N_{124}) \equiv (1,1)~mod~2,\\
    i_{c_{2}}(4) &=& (N_{234}+N_{124},N_{134}+N_{124}) \equiv (0,0)~mod~2.
    \end{eqnarray*}
    Therefore $w_{(1,2,3)}(\beta) = (1,1)_{4}(0,0)_{4} \neq 1$ in $F_{4}^{3}$ and hence the element $a_{124}a_{123}a_{124}a_{123}$ obtained in Example~\ref{exa:nontrivial-el} is non-trivial.
\end{example}


\section{Further research}\label{sec:further}
According to the result of W.P. Thurston~\cite{W-Thurston}, all knots other than satellite and torus are hyperbolic, i.e., their complement has a unique complete hyperbolic metric of constant curvature. The volume of the complement of a hyperbolic knot with respect to this metric is an invariant of the hyperbolic knot.
 
One of the very important conjectures in the knot theory is the volume conjecture relating the volume of the knot complement to some coloured Jones polynomials of this knot \cite{Kashaev,MurakamiMurakami}.
 
In the papers \cite{ManturovWan} where the photography method was introduced,
the second named author and Z. Wan noticed that volumes, areas, and other geometrical
characteristics of knots can be related to invariants constructed combinatorially,
in particular, they are related to various groups studied in \cite{FKMN} and related to braid groups. Here we mention some further directions.

\subsection{Volumes of 3-manifolds and $G_{n}^{3}$}
\begin{figure}
\centering\includegraphics[width=250pt]{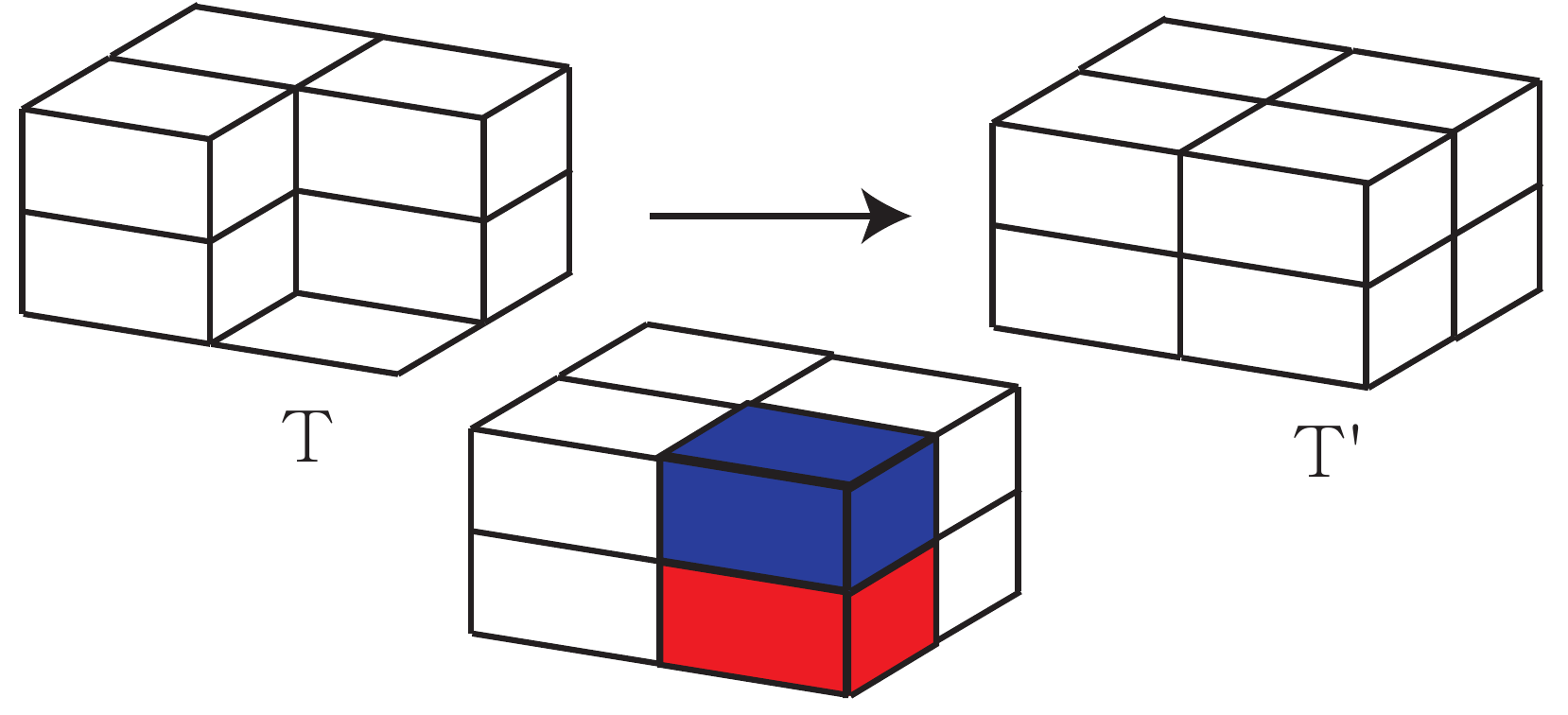}
\caption{The tiling $T'$ is obtained from $T$ by two flips. In the same time, $T'$ can be obtained from $T$ by stacking two cubes.}
\label{fig:cob-stack-cube}
\end{figure}

Let us remind that if there are two rhombile tilings on a surface and they are related by finitely many flips, then one can obtain a cobordism between them by stacking cubes, that is, one can obtain a 3-manifold with boundary components, see Fig.~\ref{fig:cob-stack-cube}. 
It follows that if a 3-manifold can be obtained by stacking cubes on a rhombile tiling, then it provides an element of $G_{n}^{3}$. In particular, the volume of the hyperbolic 3-manifold can be calculated by counting the number of hyperbolic cubes. One can expect that the number of hyperbolic cubes of the hyperbolic 3-manifold is related with the minimal length of corresponding elements in the group $G_{n}^{3}$. 

\subsection{Manturov-Nikonov indices and cubes}
In Section~\ref{sec:realisable} we construct a map from sequences of Rhombile tilings on a disc and flips to the group $G_{n}^{3}$. Since each flip is represented by a stacking of a cube, we may obtain a 3-manifold with boundaries from a sequence of rhombile tilings on a disc and flips by stacking 3-dimensional cubes. In other words, one can associated 3-manifolds obtained by stacking cubes to elements of $G_{n}^{3}$. It follows that Manturov-Nikonov indices give $\mathbb{Z}_{2}\times \mathbb{Z}_{2}$-label to cubes consisting of such 3-manifolds. In particular, since Manturov-Nikonov indices can provide a lower bound of the minimal length of the given element in $G_{n}^{3}$, one can find a lower bound of the minimal volume of 3-manifolds obtained by stacking cubes. We hope to find geometric meanings of $\mathbb{Z}_{2}\times \mathbb{Z}_{2}$-label for each cubes additionally.


\subsection{Groups $G_{n}^{k}$ and k-dimensional manifolds}
In the present paper, we studied relationship between groups $G_{n}^{3}$ and decompositions of 3-manifolds by 3-dimensional cubes. We expect that it is possible to find the relations for higher dimensional cases, that is, between $G_{n}^{k}$ and decompositions of k-manifolds by k-dimensional cubes.

\subsection{Rhombile tilings of a zonogon and braid invariants}
In \cite{Manturov-braction} it is proved that braids act on the configuration space of lines in $\mathbb{R}P^{2}$. In particular, considering Section \ref{sec:rhombi-RM3}, rhombile tilings on $\mathbb{R}P^{2}$ should correspond to elements of the configuration space of lines in $\mathbb{R}P^{2}$. Therefore, the action of a braid on $\mathbb{R}P^{2}$ might correspond to a sequence of rhombile tilings on $\mathbb{R}P^{2}$ and flips. 

On the other hand, a flip of rhombile tilings provides the equation (\ref{HS-formula}) in page \pageref{HS-formula}. That is, a sequence of rhombile tilings and flips give rise to a system of equations.

From such observations we expect that the equation (\ref{HS-formula}) should give rise to an invariant for braids.

\section*{Acknowledgements}

We are very grateful to Liliya Grunwald, Lou Kauffman and Igor Nikonov for discussion about our current work.
The first author is partially supported by the National Natural Science Foundation of China (Grant No. 12201239 and No. 12371029).

\end{document}